\documentclass{amsproc}
\usepackage{graphicx}
\font\smc=cmcsc10

\title{How to compute $\sum1/n^2$ by solving triangles}
\author{Mikael Passare}
\date{}
\address{Matematiska institutionen, Stockholms universitet\\ 
SE-10691 Stockholm, SWEDEN} \email{passare@math.su.se}

\begin{document}
    

\maketitle

\section*{Introduction}

\noindent The harmonic series $1+1/2+1/3+1/4+\ldots$ is infinitely large, whereas the series of squared terms $1+1/4+1/9+1/16+\ldots$ sums to a finite value. Indeed, the divergence of the harmonic series is easily verified by grouping its terms in packages whose lengths are increasing powers of two:
$$ 1+\frac{1}{2}+\left(\frac{1}{3}+\frac{1}{4}\right)+\left(\frac{1}{5}+\frac{1}{6}+\frac{1}{7}+\frac{1}{8}\right)+\ldots \ > \ 1+\frac{1}{2}+\frac{2}{4}+\frac{4}{8}+\ldots $$
and the boundedness of the squared series is equally simple to check:
$$1+ \frac{1}{4}+\frac{1}{9}+\frac{1}{16}+\ldots \ < \  1+ \frac{1}{1\cdot2}+\frac{1}{2\cdot3}+\frac{1}{3\cdot4}+\ldots\hskip4.2cm$$
$$\hskip4.3cm= \ 1+\left(1-\frac{1}{2}\right)+\left(\frac{1}{2}-\frac{1}{3}\right)+\left(\frac{1}{3}-\frac{1}{4}\right)+\ldots\ =\ 2\,.$$
Alternatively, one can of course also establish these well known facts by comparing with integrals:
$$\sum_{n=1}^N\frac{1}{n}\ >\ \int_0^N\frac{dx}{1+x}\ =\ \log(1+N)\ \to\ \infty\,,\ \ 
\sum_{n=1}^N\frac{1}{n^2}\ <\ 1\,+\int_1^N\frac{dx}{x^2}\ =\ 2\,-\frac{1}{N}\ \to\ 2\,.$$

The problem of computing the exact value of $\sum1/n^2$ is so classical (dating back to 1644) that it even has a name of its own: the Basel problem. The city of Basel was the hometown both of the famous brothers Jakob and Johann Bernoulli, who made serious but unsuccessful attempts to solve the problem, and to the prodigious mathematician Leonhard Euler, who found the first solution in 1735, see \cite{Eu}. Since then a great many different ways of evaluating the sum have been discovered. Fourteen proofs are collected in \cite{Chap}, and two of them have been in included in \cite{AZ}.
\smallskip

In this note we shall present a new simple geometric method for finding the exact value of $\sum_{n=1}^\infty1/n^2$ using essentially just basic trigonometry: the sine rule for triangles; combined with some elementary undergraduate analysis: area-preserving maps, see for instance \cite{Oss}.

\section*{Bipolar coordinates and the (co)sine rule}

\noindent Consider an arbitrary triangle with base of length $1$, such as the one in Fig.~1 below. The position of the top of the triangle can then be determined either by the two remaining side lengths $A$ and $B$ or, equivalently, by the two interior base angles $\alpha$ and $\beta$. One can view these data as representing two alternative sets of coordinates for the point at the top of the triangle. Thinking of the two base vertices as a pair of reference points (much like your own two eyes) it is natural to refer to $(A,B)$ and 
$(\alpha,\beta)$ as \emph{bipolar coordinates}, radial and angular respectively.

\makebox[1cm][l]{\hskip1cm\includegraphics{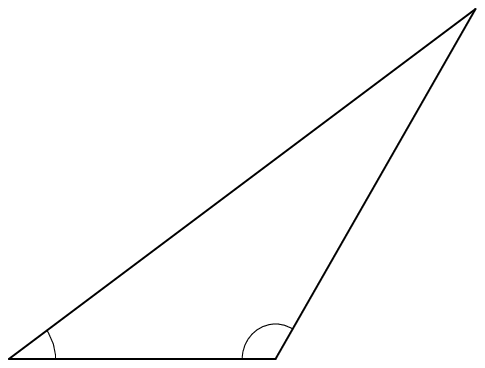}}
\begin{picture}(0,0)
\put(2.55,1.56){$\alpha$}
\put(4.2,1.85){$\beta$}
\put(3.3,1){$1$}
\put(5.5,2.5){$e^{-x}=A$}
\put(3,3.3){$B=e^{-y}$}
\put(-.18,.4){\smc Figure 1. \it The side lengths and interior angles in a triangle.}
\end{picture}\bigskip

\vskip-.2cm
Let us check what are the possible values of the two sets $(A,B)$ and $(\alpha,\beta)$ of bipolar coordinates.
For the radial coordinates $(A,B)$ the triangle inequality amounts to the obvious restriction that each side length must be less than the sum of the remaining two. That is, one has the following three inequalities: 
$$1<A+B\,,\qquad A<1+B\,,\qquad B<1+A\,.$$
These are in fact the only conditions imposed on the radial bipolar coordinates, so the collection of all possible $(A,B)$ is given by the infinite polygon $S$ depicted on the left in Fig.~2. For the angular coordinates 
$(\alpha,\beta)$, which we always measure in radians, the corresponding restriction comes from the fact that the sum of all three angles in a triangle is equal to $\pi$. The possible values of the angular coordinates $(\alpha,\beta)$ are therefore given by the half square $T$ shown on the right of Fig.~2.\medskip

\makebox[.8cm][l]{\hskip.1cm\includegraphics[width=5cm]{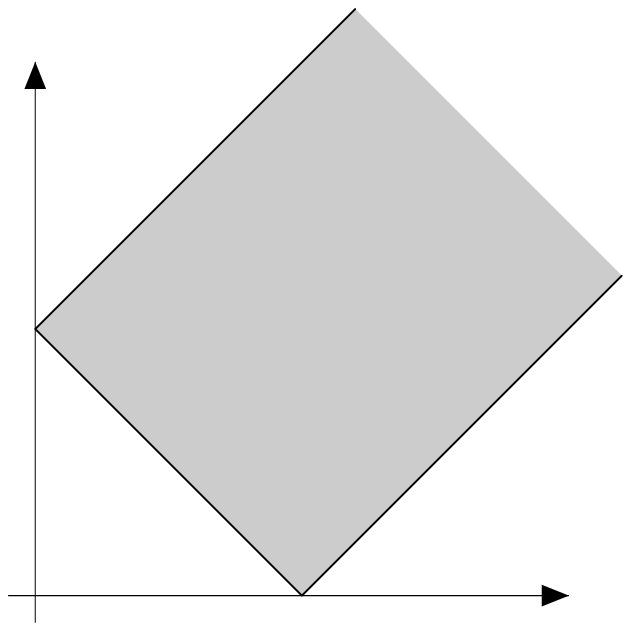}\hskip1.5cm
\includegraphics[width=5cm]{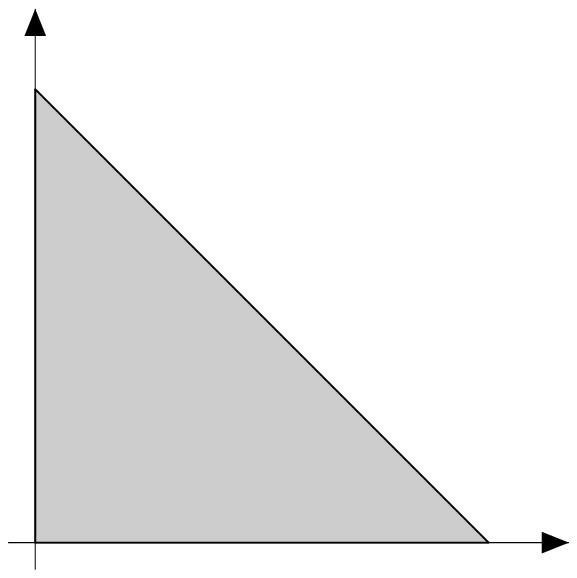}}
\begin{picture}(0,0)
\put(3.3,.55){$A$}
\put(-.4,4.05){$B$}
\put(-.25,2.45){$1$}
\put(1.62,.42){$1$}
\put(1.72,2.55){$S$}
\put(9.85,.59){$\alpha$}
\put(6.15,4.05){$\beta$}
\put(6.15,3.5){$\pi$}
\put(9.2,.48){$\pi$}
\put(7.3,1.75){$T$}
\put(-1.1,-.1){\smc Figure 2. \it The possible values for radial and angular bipolar coordinates.}
\end{picture}\bigskip

To solve a triangle means to find all its side lengths and all its interior angles, with only some of these data being known beforehand. This is precisely what one needs to do in order to pass from one set of bipolar coordinates to another. The classical tool for this is the \emph{sine rule} 
$$ \frac{\sin\alpha}{A}=\frac{\sin\beta}{B}=\frac{\sin(\pi-\alpha-\beta)}{1}=\sin(\alpha+\beta)\,,\eqno{(1)}$$
which permits one to express the side lengths $(A,B)$ in terms of the angles $(\alpha,\beta)$,
and the \emph{cosine rule}

$$A^2=1+B^2-2B\,\cos\alpha\,,\quad B^2=1+A^2-2A\,\cos\beta\,,\eqno{(2)}$$
which conversely expresses $(\alpha,\beta)$ in terms of $(A,B)$.
More precisely, the sine rule provides us with the explicit bijective map $F\colon T\to S$ given by

$$F(\alpha,\beta)=\left[\,\frac{\sin\alpha}{\sin(\alpha+\beta)}\,,\,\frac{\sin\beta}{\sin(\alpha+\beta)}\,\right],$$
whereas the cosine rule allows us to write down the inverse map 

$$F^{-1}(A,B)=\left[\,\arccos\left(\frac{1-A^2+B^2}{2B}\right),\,\arccos\left(\frac{1+A^2-B^2}{2A}\right)\,\right]\,.$$

\makebox[1cm][l]{\hskip2cm\includegraphics[width=7cm]{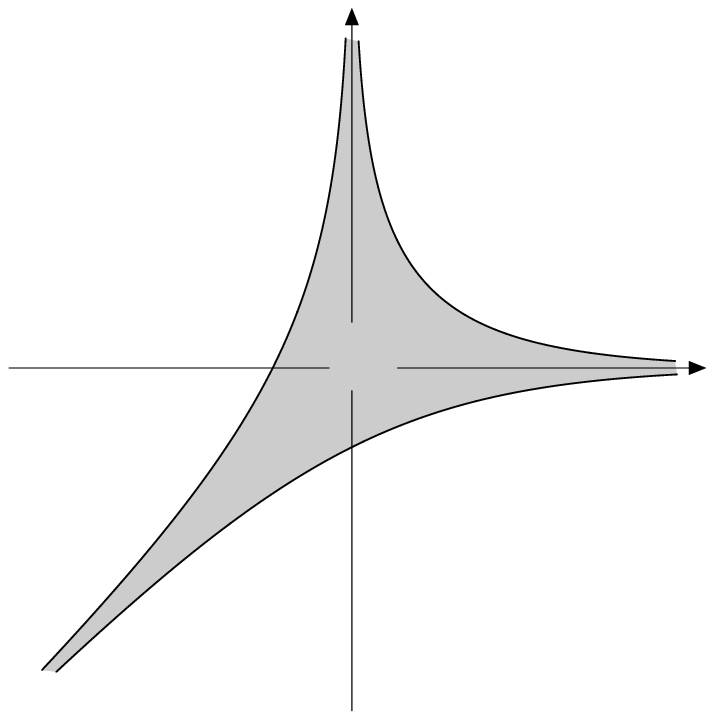}}
\begin{picture}(0,0)
\put(4.05,6.5){$y$}
\put(7.35,3.2){$x$}
\put(4.35,3.45){$U$}
\put(5.3,4.3){$e^{-x}+e^{-y}=1$}
\put(-1,-.2){\smc Figure 3. \it The set of values of the logarithmic side lengths: an amoeba.}
\end{picture}\bigskip

It will now be useful to make a logarithmic change of scale for the side lengths in Fig.~1, and thus to replace the radial bipolar coordinates $(A,B)$ by the new coordinates $(x,y)$ given by
$x=-\log A$, $y=-\log B$. The image of the set $S$ under this change of coordinates, that is, the result of drawing it on a log-log graphpaper, is the set $U$ shown in Fig.~3. This set $U$, with its ``$\,$tentacles" going off to infinity, is an example of what is nowadays known as an \emph{amoeba} in complex geometry, see \cite{PT} and \cite{Viro}. In the first quadrant its boundary is given by the curve $e^{-x}+e^{-y}=1$,
which can also be represented as $y=-\log(1-e^{-x})$. 
\smallskip

The advantage of having passed to logarithmic coordinates becomes apparent when we consider the
composed bijective map $G\colon T\to U$ given by
$$G(\alpha,\beta)=\left[\,\log\left(\frac{\sin(\alpha+\beta)}{\sin\alpha}\right),\,\log\left(\frac{\sin(\alpha+\beta)}{\sin\beta}\right)\,\right].$$
\smallskip

\noindent{\smc Theorem:} {\it The Jacobian deteminant of the map $G$ is identically equal to $1$, that is, $G$ is an area-preserving map.}
\medskip

{\it Proof:}\ Componentwise differentiation of $G$ yields the Jacobian determinant
$$\left|\begin{matrix} \cot(\alpha+\beta)-\cot\alpha & \cot(\alpha+\beta)\\
\cot(\alpha+\beta) & \cot(\alpha+\beta)-\cot\beta\\
\end{matrix}\right|=-\cot(\alpha+\beta)\bigl(\cot\alpha+\cot\beta\bigr)+\cot\alpha\cot\beta,$$
and this expression is indeed identically equal to $1$, in view of the addition formula $\cot(\alpha+\beta)=(\cot\alpha\cot\beta-1)/(\cot\alpha+\cot\beta)$.\bigskip

Since $G$ is a bijection, that is, a one-to-one map, the theorem has the following immediate and remarkable consequence.
\medskip

\noindent{\smc Corollary:} 
\hskip2cm$\text{Area}(U)=\text{Area}(T)=\pi^2/2$.
\bigskip

The threefold shape of the amoeba $U$ suggests the one should let its asymptotes divide it into three parts $U_0$, $U_1$ and $U_2$, as indicated on the left in Fig.~4. The points $(x,y)$ on either of the asymptotes correspond precisely to the side lengths $(A,B)=(e^{-x},e^{-y})$ of the isosceles triangles,
for which the interior angles $(\alpha,\beta)$ satisfy one of the conditions $\alpha=\beta$, $\alpha=\pi-\alpha-\beta$ or $\beta=\pi-\alpha-\beta$, and these are just the equations of the three medians of the triangle $T$. (In particular, the origin $(x,y)=(0,0)$ represents the equilateral triangle, with $\alpha=\beta=\pi/3$.) In other words, the map $G$ sends the medians of $T$ to the asymptotes of $U$, and one obtains a corresponding polygonal subdivision of $T$ into three parts $T_0$, $T_1$ and $T_2$, shown on the right in Fig.~4.

\makebox[.8cm][l]{\hskip-.2cm\includegraphics[width=6.5cm]{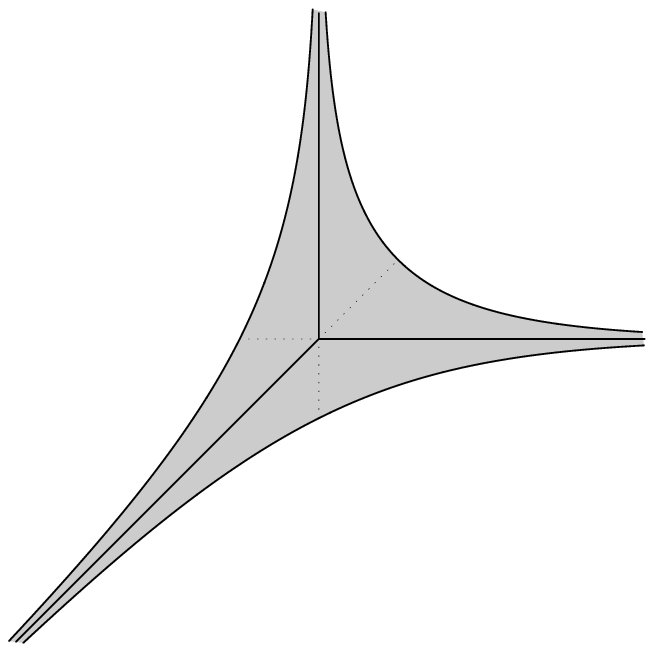}\hskip.2cm
\vbox{\includegraphics[width=5cm]{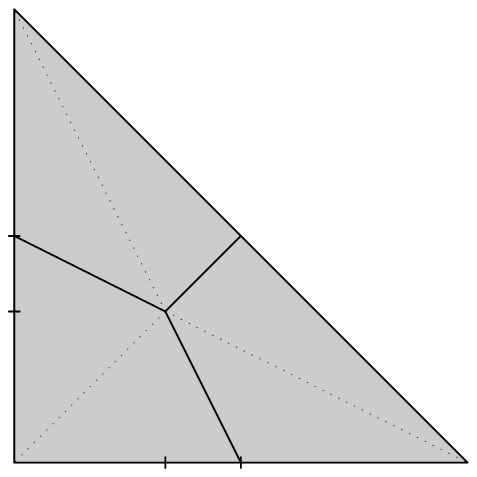}\vskip1.5cm}}
\begin{picture}(0,0)
\put(2.4,3.48){$U_0$}
\put(1.72,3.4){$U_1$}
\put(2.05,2.85){$U_2$}
\put(7.4,2.8){$T_0$}
\put(8.55,2.8){$T_1$}
\put(7.45,3.9){$T_2$}
\put(7.88,1.9){$\frac{\pi}{3}$}
\put(8.35,1.9){$\frac{\pi}{2}$}
\put(6.73,3.19){$\frac{\pi}{3}$}
\put(6.73,3.67){$\frac{\pi}{2}$}
\put(1.51,.1){\smc Figure 4. \it All six subsets have equal area.}
\end{picture}\bigskip

\noindent{\smc Theorem:} {\it The area-preserving bijective map $(x,y)\mapsto(-y,x-y)$ permutes the amoeba subsets cyclically: $U_0\mapsto U_1\mapsto U_2\mapsto U_0$.}
\medskip

{\it Proof:}\ To check that $U_0\mapsto U_1$ one can simply write down the three defining inequalities
for $U_0$ in the new coordinates $(\xi, \eta)=(-y,x-y)$:
$$x>0\ \Leftrightarrow\ \xi<\eta\,,\qquad y>0\ \Leftrightarrow\ \xi<0\,,\qquad 1<e^{-x}+e^{-y}\ \Leftrightarrow\ e^{-\xi}<e^{-\eta}+1\,,$$ 
and observe that the new inequalities are precisely the ones that define $U_1$. The rest of the theorem is proved analogously. \bigskip

It follows immediately from this theorem that $U_0$, $U_1$ and $U_2$ all have the same area, which must then be equal to $\pi^2/6$, the total area of $U$ being $\pi^2/2$. So, either using again that $G$ is area-preserving and that $G(T_0)=U_0$, or directly calculating the area of the simple polygon $T_0$, one obtains the following conclusion.
\medskip

\noindent{\smc Corollary:} 
\hskip2cm$\text{Area}(U_0)=\text{Area}(T_0)=\pi^2/6$.
\bigskip

Notice that the points in $T_0$ and $U_0$ correspond to triangles in which the longest side is the base of length $1$.
\bigskip\bigskip

\section*{Spreading and piling of squares}

\noindent Geometrically, the Basel problem amounts to computing the total area of a collection of squares with decreasing side lengths $1$, $1/2$, $1/3$, $1/4$, ..., and in order to acheive this, it will be advantageous to spread out each such square by means of an exponential function. More precisely, one observes that a square with side length $1/n$ has the same area as the part in the first quadrant lying under the exponential curve $y=e^{-nx}/n$, see Fig.~5.

\vskip-1.5cm
\makebox[1cm][l]{\hskip1cm\includegraphics[width=10cm]{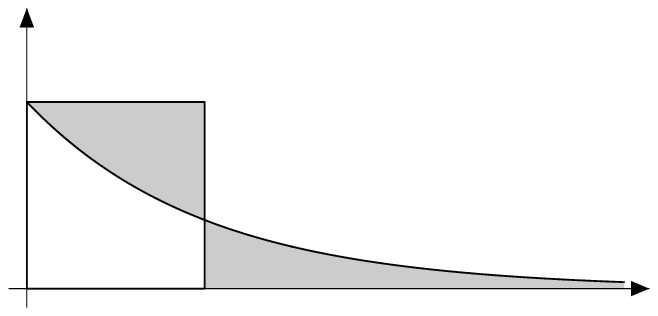}}
\begin{picture}(0,0)
\put(.6,3.4){$\frac{1}{n}$}
\put(3.12,.6){$\frac{1}{n}$}
\put(3.9,2){$y=e^{-nx}/n$}
\put(8.7,.88){$x$}
\put(.7,4.6){$y$}
\put(1.5,-.2){\smc Figure 5. \it The two areas are the same.}
\end{picture}\bigskip

\noindent Indeed, the area under the curve is readily calculated by means of the integral 
$$\frac{1}{n}\int_0^\infty e^{-nx}\,dx=\frac{1}{n}\left[-\frac{e^{-nx}}{n}\right]_0^\infty=\frac{1}{n^2}\,.$$

\noindent (Incidentally, this elementary fact can also be deduced by observing that the area-preserving bijection 
$(x,y)\mapsto((1-e^{-nx})/n,ye^{nx})$ maps the area under the curve $y=e^{-nx}/n$ to the square $0<x,y<1/n$.)
\medskip

The next step will be to pile these spread out squares on top of each other. In order to understand what such a process will yield, let us take a new look at the harmonic series. Even though the harmonic series itself is not convergent, it can be made to converge by introducing powers in the numerators.  
\bigskip

\noindent{\smc Theorem:} {\it For any positive real number $t<1$ the power series $t+t^2/2+t^3/3+t^4/4+\ldots$ is convergent and its sum is equal to $-\log(1-t)$.}
\medskip

{\it Proof:}\ Differentiation gives $$D\Big(t+\frac{t^2}{2}+\frac{t^3}{3}+\ldots+\frac{t^N}{N}\Big)=1+t+t^2+\ldots+t^{N-1}=\frac{\,\,\ 1-t^N}{1-t}\,,$$
which means that
$$t+\frac{t^2}{2}+\frac{t^3}{3}+\ldots+\frac{t^N}{N}=\int_0^t\frac{\,\,\ 1-\tau^N}{1-\tau}\,d\tau=\int_0^t\frac{d\tau}{1-\tau}-\int_0^t\frac{\tau^N\,d\tau}{1-\tau}\,.$$
Here the first integral on the right has the desired value $-\log(1-t)$ while the second integral is less than
$$t^N\int_0^t\frac{d\tau}{1-\tau}=-t^N\log(1-t)\,,$$ which tends to zero as $N\to\infty$. \bigskip

Writing $t=e^{-x}$, so that $t^n=e^{-nx}$, we thus see that for every $x>0$ the series $\sum e^{-nx}/n$ has the sum $-\log(1-e^{-x})$. But piling the spread out squares on top of each other, as indicated on the right in Fig.~6, precisely amounts to summing the functions $e^{-nx}/n$, so from the theorem one deduces the following fact.
\medskip

\noindent{\smc Corollary:} {\it The infinite pile of spread out squares exactly covers the set $U_0$ defined by the inequalities $x>0$, $y>0$, $y<-\log(1-e^{-x})$.}

\vskip-.5cm
\makebox[1cm][l]{\hskip1cm\includegraphics[width=10cm]{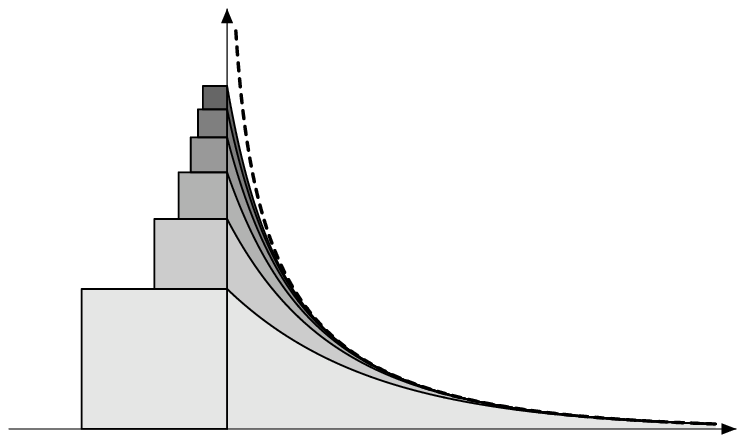}}
\begin{picture}(0,0)
\put(3.5,3.5){$e^{-x}+e^{-y}=1$}
\put(8.9,.61){$x$}
\put(2.34,6.05){$y$}
\put(-1,0){\smc Figure 6. \it Building two infinite piles. Equally shaded areas are the same.}
\end{picture}\bigskip

Since the convergence is not uniform in $x$ a little argument may be in order here: One easily verifies that for any compact subset $K$ of $U_0$ there is some finite pile $0<y<e^{-x}+e^{-2x}/2+\ldots+e^{-Nx}/N$ which contains $K$. On the other hand, it is clear that every finite pile is contained in $U_0$.
\bigskip

\noindent{\smc Conclusion:} {\it The series $1+1/4+1/9+\ldots$ of squared reciprocals of positive intergers has the exact sum $\pi^2/6$.}
\vfill\eject

In order to emphasize the conciseness of the solution to the Basel problem that we have just presented, let us recapitulate it again in formulas:

$$ \sum_{n=1}^\infty\frac{1}{n^2}= \sum_{n=1}^\infty\,\,\int_0^\infty\frac{e^{-nx}}{n}\,dx=-\int_0^\infty\!\!\log(1-e^{-x})\,dx=\int_{U_0}\!\!dx\,dy=\int_{T_0}\!\!d\alpha\,d\beta=\frac{\pi^2}{6}\,.$$
\smallskip

\noindent Here the notations $U_0$ and $T_0$ again refer to the sets appearing earlier, see  Fig.~4.
\bigskip\bigskip

\section*{Epilogue: the complex logarithm explains the mystery}

\noindent The most surprising part of the arguments in the preceding sections is probably the fact that the map $G$, connecting the logarithmic side lengths $(x,y)$ and the angles $(\alpha,\beta)$, turned out to be area-preserving. Why on earth should the Jacobian of $G$ be identically equal to $1\,$? 
\smallskip

In order to shed some light on this enigma, we shall take a look at logarithms of complex numbers. Vaguely put, it is their marvelous property 
$$\log z=\log|z|+i\arg z\,,$$
combining (the logarithm of) the modulus $|z|$ with the argument $\arg z$, that explains the presence of a simple relation between the logarithmic side lengths $(x,y)$ and the angles $(\alpha,\beta)$. 
\smallskip

Let us now be more precise, and
consider the linear equation $1+z+w=0$, with the unknowns $z$ and $w$ being complex numbers. Depicting $1$, $z$ and $w$ as vectors in the plane, one can interpret the fact that they sum to zero as saying the three vectors should form a closed triangle, quite similar to the one in Fig.~1.

\makebox[1cm][l]{\hskip1cm\includegraphics{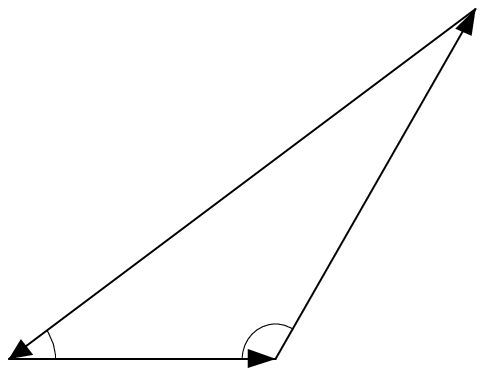}}
\begin{picture}(0,0)
\put(2.55,1.56){$\alpha$}
\put(4.2,1.85){$\beta$}
\put(3.3,1){$1$}
\put(5.5,2.5){$|z|=A$}
\put(3,3.3){$B=|w|$}
\put(-.18,.4){\smc Figure 7. \it The triangle formed by the vectors $1$, $z$ and $w$.}
\end{picture}\bigskip

We now get $1+z+w=1+Ae^{i(\pi-\beta)}+Be^{i(\alpha+\pi)}=0$, by choosing the arguments $\arg z$ and $\arg w$ to lie between $0$ and $2\pi$, and taking real and imaginary parts of this complex equation, we see that it amounts to the two real equations
\begin{equation*}
\hskip3.7cm
\begin{cases}
A\sin\beta-B\sin\alpha &=\ 0\ ,  \hskip3.7cm (3) \\
1-A\cos\beta-B\cos\alpha &=\ 0\ . \hskip3.75cm (4) \\
\end{cases}
\end{equation*}
It is an amusing exercise to verify that this system of equations is in fact equivalent to the sine and cosine rules (1) and (2). For instance, to deduce the cosine rule from (3) and (4), one can first re-write (4) as $A\cos\beta=1-B\cos\alpha$. Squaring both sides, and using (3) to replace the squared left hand side $A^2\cos^2\beta=A^2(1-\sin^2\beta)$ by $A^2-B^2\sin^2\alpha$, one obtains the cosine rule $A^2=1+B^2-2B\cos\alpha$.
\smallskip

Letting the arguments $\arg z$ and $\arg w$ be more general, it is natural to ``lift" the complex line $1+z+w=0$ to the exponential complex curve
$$X=\bigl\{\,(s,t)\in\mathbf{C}^2\,;\ 1+e^s+e^t=0\,\bigr\}\,.$$
Denote by $P$ and $Q$ respectively the restrictions to $X$ of the two linear projections 
$$ (s,t)\longmapsto (\mathrm{Re}\,s,\mathrm{Re}\,t)\quad \text{and} \quad (s,t)\longmapsto (\mathrm{Im}\,s,\mathrm{Im}\,t)$$
on the real and imaginary parts, and let $\widetilde U=P(X)$, $\widetilde T=Q(X)$ be the corresponding images in $\mathbf{R}^2$. Observe that $Q$ is invertible on the interior of 
$\widetilde T$, so we can introduce a map $\widetilde G=P\circ Q^{-1}$ going from imaginary parts to real parts. This map $\widetilde G$ is closely related to our previous map $G$. In fact, one has the identity
$$G=\left[
\begin{matrix} 
-1 & \,0\, \\
\,\,\,0 & \!\!-1\, \\
\end{matrix}\right]
\circ\widetilde G\circ
\left[
\begin{matrix} 
\,\,0 & \!\!-1\, \\
\,\,1 & \,\,0 \, \\
\end{matrix}\right]\,.$$
Here the two linear maps (represented by their matrices) are both area-preserving (since the determinants equal $1$), so what remains to be explained is why $\widetilde G$ has its Jacobian identically equal to $1$.
\smallskip

Fix any point $p\in X$ and let $T_p$ denote the (real two-dimensional) tangent space to $X$ at $p$. Since $P$ and $Q$ are (restrictions to $X$ of) linear maps, their derivatives $P'(p)$ and $Q'(p)$ coincide with (restrictions to $T_p$ of) the maps themselves. The Jacobian of $\widetilde G$ is therefore equal to the determinant of $P\circ Q^{-1}$, where we now use the same notation $P$, $Q$ to mean the real and imaginary projections $T_p\to\mathbf{R}^2$. 
\smallskip

Notice next that $T_p$ is in fact a (one-dimensional) complex subspace of $\mathbf{C}^2$, so coordinate-wise multiplication by $i$ defines a linear map $J$ on $T_p$. Since $J$ just amounts to a rotation (by an angle $\pi/2$) around the origin it is clear that $J$ is area-preserving and thus has determinant $1$. Moreover, one has $J=Q^{-1}\circ P$, so $P\circ Q^{-1}=
Q\circ(Q^{-1}\circ P)\circ Q^{-1}$ also has determinant $1$, which means that $\widetilde G$ has Jacobian equal to $1$.

\end{document}